\input{STU.STY}
\usepackage{xcolor}
\newfont{\cyr}{wncyr8}
\newfont{\cyb}{wncyr8}
\newtheorem{thm}{Theorem}[section]
\newcommand{\bthm}{\begin{thm}}
\newcommand{\ethm}{\end{thm}}

\newtheorem{prop}[thm]{Proposition}
\newcommand{\bprp}{\begin{prop}}
\newcommand{\eprp}{\end{prop}}

\newtheorem{fact}[thm]{Fact}
\newcommand{\bfct}{\begin{fact}}
\newcommand{\efct}{\end{fact}}

\newtheorem{prob}[thm]{Problem}
\newcommand{\bprb}{\begin{prob}}
\newcommand{\eprb}{\end{prob}}

\newtheorem{quest}[thm]{Question}
\newcommand{\bqtn}{\begin{quest}}
\newcommand{\eqtn}{\end{quest}}

\newtheorem{lem}[thm]{Lemma}
\newcommand{\blem}{\begin{lem}}
\newcommand{\elem}{\end{lem}}

\newtheorem{claim}[thm]{Claim}
\newcommand{\bclm}{\begin{claim}}
\newcommand{\eclm}{\end{claim}}

\newtheorem{cor}[thm]{Corollary}
\newcommand{\bcor}{\begin{cor}}
\newcommand{\ecor}{\end{cor}}

\newtheorem{conj}[thm]{Conjecture}
\newcommand{\bcnj}{\begin{conj}}
\newcommand{\ecnj}{\end{conj}}

\theoremstyle{definition}
\newtheorem{defn}[thm]{Definition}
\newcommand{\bdfn}{\begin{defn}}
\newcommand{\edfn}{\end{defn}}

\newtheorem{spec}[thm]{Specializing}
\newcommand{\bspc}{\begin{spec}}
\newcommand{\espc}{\end{spec}}

\theoremstyle{remark}
\newtheorem{rem}[thm]{Remark}
\newcommand{\brem}{\begin{rem}}
\newcommand{\erem}{\end{rem}}

\newtheorem{cnv}[thm]{Convention}
\newcommand{\bcnv}{\begin{cnv}}
\newcommand{\ecnv}{\end{cnv}}

\theoremstyle{Example}
\newtheorem{exam}[thm]{{\bf Example}}
\newcommand{\bexm}{\begin{exam}}
\newcommand{\eexm}{\end{exam}}

\newtheorem{exercise}[thm]{Exercise}
\newcommand{\bexr}{\begin{exercise}}
\newcommand{\eexr}{\end{exercise}}

\newtheorem{thmy}{\textbf{Theorem}}
\newenvironment{thmx}{\stepcounter{thm}\begin{thmy}}{\end{thmy}}

\usepackage{pdfpages}
\usepackage{amsxtra}

\begin{document}
\label{begin-art}

\title{Duality and Baire like properties in topological abelian groups}

\date{\today}

\author[M. Ferrer, S. Hern\'andez, I. Sep\'ulveda]
{M. Ferrer, S. Hern\'andez, I. Sep\'ulveda}
\address{Universitat Jaume I, Departamento de Matem\'{a}ticas,
Campus de Riu Sec, 12071 Castell\'{o}n, Spain.}
\email{mferrer@mat.uji.es, hernande@mat.uji.es, isepulve@mat.uji.es}

\thanks{We acknowledge partial financial support by the Institute of Mathematics and its Applications of Castell\'o
(Universitat Jaume I)}

\begin{abstract}
This article investigates the duality of abelian topological groups by delving
into the interplay between the algebraic and topological properties of a group
$G$ and those of its dual group. In particular, extending the notion of a
$g$-barrelled group, we introduce and study the classes of $\kappa$-barrelled
and $\omega$-barrelled groups. Our first main result establishes that every
$\omega$-barrelled group is $\aleph_0$-barrelled (i.e., if every weakly
convergent sequence in the dual group is equicontinuous, then every metrizable
weakly compact subset of the dual group is also equicontinuous). From this, we
deduce that the free topological abelian group $A(K)$ over a metrizable compact
space $K$ is determined by the subgroup generated by any of its dense subsets.

Furthermore, we prove that if $G$ is a $\sigma$-compact, $\omega$-barrelled
group, then its dual group $\widehat{G}_\kappa$ equipped with the compact-open
topology is sequentially complete (complete if $G$ is hemicompact).
We also solve an open question posed by Trigos-Arrieta by constructing an
explicit example of an $\omega$-barrelled metrizable group that fails to be
$g$-barrelled. Finally, we investigate the dualization of groups lacking
infinite compact subsets, proving that for a totally bounded abelian group $G$,
every compact subset of $G$ is finite if and only if its dual group $\widehat{G}$,
endowed with the finite-open topology, is unordered Baire-like.
\end{abstract}


\maketitle \setlength{\baselineskip}{24pt}

\section{Introduction}

The classical Pontryagin--van Kampen duality theorem establishes a full duality for locally compact abelian (LCA) groups, yielding a powerful framework where topological and algebraic properties closely reflect one another (see, e.g., \cite{ausdikgio, HR, hofmor_compactgroups, morr, mackey, rudin}). Beyond the LCA setting, however, this tight correspondence often breaks down. Accordingly, a main current line of research in the theory of topological abelian groups aims to understand how algebraic and topological properties interact under duality in more general
classes of groups.

In this paper, we investigate this interplay by focusing on the extension of the uniform boundedness theorem to general topological abelian groups. Specifically, our primary objective is to classify those groups for which certain compact subsets of their duals are equicontinuous.

While Baire's category property solves this problem for certain classes of groups, it fails to do so in general. To address this gap, Chasco, Martín-Peinador, and Tarieladze \cite{ChMPVa:99} introduced the notion of a \emph{$g$-barrelled group}, defined as a group in which every compact subset of $\widehat{G}_p$ (i.e., the dual group $\wG$ equipped with the topology of pointwise convergence) is equicontinuous. Following a previous strategy of Saxon \cite{Saxon:MA1972} for topological vector spaces, we investigate this property by utilizing variants of Baire's category property.

Further motivation comes from the study of totally bounded groups. In \cite{ChascoDominguezTkachenko,AusDik:2021,FerHerTka}, the authors studied groups $G$ whose dual $\widehat{G}_p$ contain no infinite compact subsets, leading to the conjecture that if $G$ has the Baire property, then $\widehat{G}_p$ contains no infinite compact subsets---a result subsequently confirmed in \cite{FerHerSepTri}. Conversely, the precise dual characterization of groups lacking infinite compact subsets remains an open problem. One of our main goals is to clarify this and related structural questions.

To this end, we introduce and study several natural generalizations of $g$-barrelledness. From now on, the \emph{weight} of a topological space $X$, denoted $\we(X)$, is the minimal cardinality of a basis of open sets for the topology of $X$.
We say that a group $G$ is \emph{$\kappa$-barrelled} if every compact subset $K \subseteq \widehat{G}_p$ with $\mathrm{w}(K) \leq \kappa$ is equicontinuous. In particular, the group $G$ is called \emph{$\omega$-barrelled} (or \emph{sequentially $g$-barrelled}) if every convergent sequence in $\widehat{G}_p$ is equicontinuous. Along similar lines, a subset $U \subseteq G$ is said to be \emph{$\kappa$-dense} if there exists a subset $A \subseteq G$ with $|A| \leq \kappa$ such that $G = A + U$. In particular, a topological Abelian group is \emph{$\omega$-narrow}
if every neigborhood of $0_G$ is $\aleph_0$-dense.

Let $\mathcal{K}(\widehat{G}_p)$ denote the family of all compact subsets of $\widehat{G}_p$. A quasi-convex subset $B \subseteq G$ is called a \emph{$\mathcal{K}$-neighborhood} of $0_G$ if there exists a compact set $K \in \mathcal{K}(\widehat{G}_p)$ such that $B = K^\triangleleft$, where $K^\triangleleft$ denotes the polar of $K$. We say that $G$ is \emph{unordered Baire-like} if for every sequence $\{B_n\}_{n<\omega}$ of $\mathcal{K}$-neighborhoods of $0_G$ whose union is $G$,
at least one $B_n$ has a non-empty interior. We now present our main findings below.\medskip

Looking at the definitions listed above, one might conjecture that every $\go$-barrelled metrizable group must be $g$-barrelled (see \cite{Trigos:2025}, where this question was explicitly proposed). In Section \ref{Javier's question}, we provide a negative answer to this question by constructing a counterexample that disproves this conjecture. However, the following fact holds.

\begin{thmx}\label{Thm_A} (Thm-A)
Every $\omega$-barrelled group is $\aleph_0$-barrelled.
\end{thmx}
\mkp

\brem\label{Rem/28/4} (Rem/28/4) According to Comfort, Raczkowski, and Trigos-Arrieta \cite{comractri:04}, a group $G$ is said to be
\emph{determined by a subgroup} $H$ if the dual groups of $G$ and $H$, equipped with the compact-open topology, are topologically isomorphic. The group $G$ is said to be \emph{determined} if it is determined by every dense subgroup. The cornerstone of this topic was
established by Aussenhofer \cite{Aussenhofer} and, independently, Chasco \cite{chasco}, who proved that every metrizable abelian group $G$ is determined. Subsequently, it was proved in \cite{FHU:jmaa} that a compact group (abelian or not) is determined if and only if
it is metrizable. The following Corollary is a direct consequence of Theorem \ref{Thm_A}. It provides new examples of
non-metrizable groups that are determined by dense subgroups.
\erem
\bkp

\noindent {\bf Corollary A}\label{Thm_B} (Thm-B)
Let $K$ be a metric compact space and let $D$ be a dense subset of $K$. Then the free topological abelian group $A(K)$
is determined by the subgroup generated by $D$.
\mkp

A variation of the proof of Theorem \ref{Thm_A} allows us to establish the following fact.

\begin{thmx}\label{Thm_C} (Thm-C)
Let $G$ be a $\sigma$-compact, $\go$-barrelled group. Then its dual group $\widehat{G}_\sK$, equipped with the compact-open topology, is sequentially complete. Therefore, if $G$ is hemicompact then $\widehat{G}_\sK$ is complete.
\end{thmx}
\mkp

Another question that remains open in the literature is the characterization of totally bounded abelian groups containing no infinite compact subsets. The following result generalizing \cite[2.12]{FerHerSepTri}, addresses this question.

\begin{thmx}\label{Thm_D} (Thm-D)
Let $G$ be a totally bounded abelian group. Then every compact subset of $G$ is finite if and only if its dual group $\widehat{G}_p$, equipped with the finite-open topology, is unordered Baire-like.
\end{thmx}
\mkp

The paper is organized as follows: In Section 2, we introduce the algebraic and topological tools that will be used throughout the article. In Section 3, we develop several technical results and establish relevant preliminary implications. In Section 4, we study Baire-like groups and examine their basic properties. Section 5 contains the main results of the paper, where we focus on $\omega$-barrelled groups and analyze some of their remarkable properties and structural equivalences. Section \ref{Javier's question} is devoted to providing an explicit example of an $\omega$-barrelled metrizable group that is not $g$-barrelled. Finally, in Section 6, we study unordered Baire-like groups, leading to a characterization of $g$-barrelled groups within the class of totally bounded abelian groups.

\section{Basic facts}

Unless stated otherwise, all groups are assumed to be abelian and all topological spaces are Tychonoff. We denote the group of integers by $\ZZ$; depending on the context, $\ZZ$ may be equipped with the discrete topology or viewed purely algebraically as an abstract group. Similarly, $\ZZ_q$ denotes the cyclic group of order $q \in \NN$. Set $\omega := \{0, 1, 2, \dots\}$. If $S$ is a set, ${\rm Fin}(S)$ denotes the family of all finite subsets of $S$.

Our model for the torus $\TT$ is the interval $(-1/2, 1/2]$ equipped with the addition modulo $1$ and the topology inherited from $\RR$ with the endpoints $-1/2$ and $1/2$ identified, unless we explicitly regard it as a topological group in another standard representation.

Let $G$ be an abelian group. For a subset $A \subseteq G$ and a natural number $n \in \NN$, we define the set of at most $n$-fold sums from $A$ by
\[
(n)A := \{a_0 + \dots + a_j : 0 \leq j \leq n, \ a_1, \dots, a_j \in A\},
\]
where $j=0$ corresponds to the identity element $0_G$. The subgroup generated by $A$, denoted by $\langle A \rangle$, is the smallest subgroup of $G$ containing $A$; if $A$ is a singleton $\{a\}$, we simply write $\langle a \rangle$. Finally, $\mathrm{Tor}(G)$ denotes the torsion subgroup of $G$.

\subsection{Pontryagin duality}

Let $(G, \tau)$ be a topological abelian group with underlying group $G$ and topology $\tau$. A \emph{character} on $(G, \tau)$ is a continuous homomorphism from $G$ into the torus $\TT$. The pointwise sum of two characters is again a character; thus, the set $\wG$ of all characters forms an abelian group under pointwise addition. When $\wG$ is equipped with the compact-open topology, it becomes a topological group $(\wG, \widehat{\tau})$, known as the \emph{dual group} of $(G, \tau)$.

There is a natural evaluation homomorphism $E_G \colon G \longrightarrow \widehat{\wG}$ from $G$ into its bidual group, which is not necessarily continuous, injective, or surjective. Evidently, the evaluation map $E_G$ is injective if and only if $G$ is maximally almost periodic (MAP). A topological abelian group $(G, \tau)$ satisfies \emph{Pontryagin duality} or is \emph{Pontryagin reflexive} (\emph{P-reflexive} for short) if the evaluation map $E_G$ is a topological isomorphism onto $\widehat{\wG}$. If $E_G$ is merely a group isomorphism (i.e., bijective but not necessarily a homeomorphism), we say that $G$ is \emph{P-semireflexive}.

\subsection{Topologies on the dual group induced by bornologies}

In this paper, a \emph{bornology} $\mathcal{B}$ on an abelian topological group $G$
designates an ideal of subsets of $G$ satisfying the following properties:
\begin{itemize}
\item[(i)] $\mathrm{Fin}(G) \subseteq \mathcal{B}$;
\item[(ii)] $A \in \mathcal{B}$ and $B \subseteq A$ imply $B \in \mathcal{B}$;
\item[(iii)] $A, B \in \mathcal{B}$ imply $A \cup B \in \mathcal{B}$;
\item[(iv)] $A, B \in \mathcal{B}$ imply $A + B \in \mathcal{B}$.
\end{itemize}
Here, the symbol $\mathrm{Fin}(G)$ denotes the collection of all finite subsets of $G$.
On the other hand, statements (ii) and (iii) are already contained in the definition of an ideal.

\bigskip

Let $\sB(G)$ be a bornology on $G$. For a subset $A \subseteq G$ and a positive real number $\epsilon > 0$, we define
\[
[A, \epsilon] = \{\chi \in \wG : |\chi(a)| \leq \epsilon \text{ for all } a \in A\}.
\]

It is a well-known fact that the family of sets $[K, \epsilon] \subseteq \wG$, where $K \in \sB(G)$ and $\epsilon > 0$, forms a neighborhood base at the trivial character, thereby defining a group topology on $\wG$ which we denote by $t(\sB)$. We write $\wG_\sB$ for the topological group obtained in this manner.

For example, if we consider the bornology $\sK(G)$ generated by all compact subsets of $G$, we obtain the classical Pontryagin-van Kampen \emph{dual} group $\wG_{\sK(G)}$ (or simply $\wG_\sK$). Taking the compact bornology $\sK(\wG_\sK)$ on this dual group, we obtain a new topological group, namely the \emph{bidual} group $\widehat{(\wG_\sK)}_{\sK(\wG_\sK)}$ (denoted by $\widehat{\wG}_\sK$ for short).
Within this terminology, the group $G$ is said to be (Pontryagin-van Kampen) \emph{reflexive} if the evaluation map
\[
E_G : G \to \widehat{\wG}_\sK
\]
is a topological isomorphism. 

Let $K \in \sB$. For each $n \in \NN$, and using (iv) in the bornology definition, we have $[(n)K, 1/4] \subseteq [K, 1/4n]$. Thus, the collection of sets $[K, 1/4]$ for $K \in \sB$ also forms a neighborhood base for $\wG_\sB$ at the trivial character.

\brem
Observe that if $G$ is a topological abelian group and $\sB$ is a bornology on the dual group $\wG$, then we can replicate the notions introduced above to define the topology $t(\sB)$ on $G$. Hence, given a subset $X \subseteq \wG$ and a positive real number $\epsilon > 0$, we define
\[
[X, \epsilon] = \{g \in G : |\chi(g)| \leq \epsilon \text{ for all } \chi \in X\}.
\]
Again, the sets $[K, \epsilon] \subseteq G$ for $K \in \sB(\wG)$ and $\epsilon > 0$ form a neighborhood base at the identity element, defining a group topology on $G$.
\erem

\bdfn
For a subset $A \subseteq G$, we define its polar by $A^\rhd = [A, 1/4]$. Similarly, for a subset $X \subseteq \wG$, its polar (in $G$) is defined by $X^\lhd = \{g \in G : |\chi(g)| \leq 1/4 \text{ for all } \chi \in X\}$.
\edfn

\blem[{\cite[Lemma 2.2]{Noble:tams}}]\label{Upol}
For each neighborhood $U$ of $0$ in $G$, we have $U^\rhd \in \sK(\wG_\sK)$.
\elem

\bdfn[Vilenkin \cite{Vilenkin}]
A subset $A \subseteq G$ is called \emph{quasiconvex} if $A^{\rhd\lhd} = A$. The group $G$ is called \emph{locally quasiconvex} if it possesses a neighborhood base at the identity consisting of quasiconvex sets.
\edfn

For each $A \subseteq G$, its polar $A^\rhd$ is a quasiconvex subset of $\wG$. Thus, given a bornology $\sB$ on $G$, the topological group $\wG_{\sB}$ is locally quasiconvex for any topological group $G$ and any bornology $\sB$ on it. Moreover, because local quasiconvexity is hereditary for arbitrary subgroups, the group $G_\sB$ is locally quasiconvex for any bornology $\sB$ on $\wG$. The subset $A^{\rhd\lhd}$ is the smallest closed quasiconvex subset of $G$ containing $A$.

In the case where $G$ is a topological vector space, $G$ is locally quasiconvex in the present sense if and only if it is a locally convex topological vector space in the ordinary sense \cite{Banaszczyk}.

If $G$ is locally quasiconvex, its characters separate the points of $G$, implying that the evaluation map $E_G : G \to \widehat{\wG}_\sK$ is injective. Furthermore, for each quasiconvex neighborhood $U$ of $0$ in $G$, the polar $U^\rhd$ is a compact subset of $\wG_\sK$ (by Lemma \ref{Upol}), which means that $U^{\rhd\rhd}$ is a neighborhood of $0$ in $\widehat{\wG}_\sK$. Since $E_G[G] \cap U^{\rhd\rhd} = E_G[U^{\rhd\lhd}] = E_G[U]$, the evaluation map $E_G$ is open onto its image \cite[Lemma 14.3]{Banaszczyk}.

\bigskip
The next proposition summarizes some of these statements.

\bprp\label{pr_top}
If $\sB$ is a bornology on a topological abelian group $G$ (resp. $\wG$), then $\wG_\sB$ is a Hausdorff, locally quasiconvex group (resp. $G_\sB$ is a locally quasiconvex group).
\eprp

\subsection{Totally bounded groups}

For a topological group $G$, we denote by $\mathrm{Hom}(G, \TT)$ the set of all group homomorphisms $\phi \colon G \to \TT$, which forms a group under pointwise sum. Thus, the dual group $\wG$ consists precisely of the continuous homomorphisms within $\mathrm{Hom}(G, \TT)$. When $G$ is discrete, its dual group coincides with $\mathrm{Hom}(G, \TT)$.

A topological group $G$ is called \emph{precompact} if for every neighborhood $U$ of $0$, there exists a finite subset $F \subseteq G$ such that $G = F + U$. If, in addition to being precompact, $G$ is Hausdorff, we say that $G$ is \emph{totally bounded}. It is a well-known fact that for any topological group $G$, the finite bornology $\mathrm{Fin}(G)$ defines a totally bounded topology on $\wG$, denoted by $\wG_{\mathrm{Fin}(G)}$ (or simply $\wG_p$).

If $S$ is a subgroup of $\mathrm{Hom}(G, \TT)$, we denote by $G_{\mathrm{Fin}(S)}$ (or $G_S$ for short) the topological group obtained by equipping $G$ with the topology $\tau_S$ generated by the bornology $\mathrm{Fin}(S)$. This topology is the coarsest group topology on $G$ rendering the elements of $S$ continuous. It follows that $G_S$ is a precompact group, which is Hausdorff if and only if $S$ separates the points of $G$ (i.e., if $S$ is point-separating). This is equivalent to asserting that $S$ is dense in $\mathrm{Hom}(G, \TT)$, where the latter is compact when viewed as the dual group $\wG_p$ of a discrete group $G$. By the Comfort-Ross Theorem \cite{ComfortRoss1964}, we have $(\widehat{G_S})_p = S$.

A topological group $G$ is said to be \emph{maximally almost periodic} (MAP for short) if $\wG$ separates the points of $G$. In this case, $\wG$ is dense in the Bohr compactification $b\wG_p$, whose underlying algebraic group is precisely $\mathrm{Hom}(G, \TT)$ (see \cite[Theorem 1.12]{ComfortRoss1964}).

\section{Baire like groups}\label{S3}
In this section, we introduce and study the basic properties and relationships of several classes of groups that generalize the notion of Baire groups.
\medskip

\bdfn\label{def1}
Let $G$ be a topological abelian group and let $\sB$ be a bornology on $\wG$. A quasi-convex subset $B$ of $G$ is said to be a \emph{$\sB$-neighborhood} of $0_G$ if there is an element $K \in \sB$ such that $B = K^\triangleleft$. A \emph{nest-generated $\sB$-barrel set} (\emph{ng-$\sB$-barrel set} for short) is a sequence $\{B_n\}_{n \in \NN}$ of $\sB$-neighborhoods of $0_G$ satisfying:
\begin{enumerate}
    \item[(a)] $\bigcup_{n \in \NN} B_n$ contains an open subgroup of $G$.
    \item[(b)] $B_n \subseteq B_{n+1}$ for all $n \in \NN$.
\end{enumerate}

To simplify the notation, when $\sB = \sK(\widehat G)$, we simply say \emph{nest-generated barrel set} (\emph{ng-barrel set} for short). Extending the definition given by Saxon in \cite{Saxon:MA1972}, the group $G$ is said to be \emph{$\sB$-Baire-like} (resp., \emph{Baire-like}) if every ng-$\sB$-barrel set $\{B_n\}$ contains an element with a non-empty interior (resp., every ng-barrel set $\{B_n\}$ contains an element with a non-empty interior). Note that if $\sB$ contains all convergent sequences in $\wG_p$, then $G$ is $\omega$-barrelled (see Theorem \ref{prp1.1} below).

A group $G$ is said to be \emph{strongly $\sB$-Baire-like} (resp., \emph{strongly Baire-like}) if every ng-$\sB$-barrel set $\{B_n\}$ stabilizes; that is to say, there exists $n_0 \in \NN$ such that $B_n = B_{n_0}$ for all $n \geq n_0$ (resp., every ng-barrel set $\{B_n\}$ stabilizes).

A subset $A$ of a topological space $X$ is said to be \emph{$G_\delta$-dense} in $X$ if every non-empty $G_\delta$ subset $U$ of $X$ intersects $A$. The group $G$ is said to be \emph{dually Baire-like} if for every sequence $\{F_n\}$ of $\sK$-neighborhoods of $0_G$ such that $\bigcup_{n \in \NN} F_n$ is $G_\delta$-dense in $G$, there exists $n_0 \in \NN$ such that $(F_1 \cup \dots \cup F_{n_0})^\rhd = \{0\}$.
\edfn
\medskip

We will now examine some basic relationships between the concepts introduced above. First, we establish the following folklore lemma, the proof of which is included for the sake of completeness.
\medskip

\blem\label{lem1}
Let $\chi \in \widehat{G}$ and $g_1, g_2 \in G$ such that $|\chi(g_i)| \leq \frac{1}{8}$ for $i = 1, 2$. Then $\chi(g_1 + g_2) = \chi(g_1) + \chi(g_2)$.
\elem
\pf Since $\chi(g_1 + g_2) \in (-\frac{1}{2}, \frac{1}{2}]$ and $\chi(g_i) \in [-\frac{1}{8}, \frac{1}{8}]$ for $i = 1, 2$, it follows that $\chi(g_1) + \chi(g_2) \in [-\frac{1}{4}, \frac{1}{4}]$. Thus, we have
\[
\chi(g_1 + g_2) - \chi(g_1) - \chi(g_2) \in \left( \left( -\frac{1}{2}, \frac{1}{2} \right] + \left[ -\frac{1}{4}, \frac{1}{4} \right] \right) \cap \ZZ \subseteq \left[ -\frac{3}{4}, \frac{3}{4} \right] \cap \ZZ = \{0\}.
\]
\epf

\bthm\label{prp1.1} (prp1.1)
Let $G$ be a topological abelian group and let $\sB$ be a bornology on $\wG$. The following two assertions hold true:
\begin{enumerate}
  \item[(1)] If for every open subgroup $H$ of $G$, the dual group $\widehat{G}$ contains no non-trivial pointwise convergent sequences on $H$, then every ng-$\sB$-barrel set stabilizes and, as a consequence, the group $G$ is strongly $\sB$-Baire-like.
  \item[(2)] If $G$ is $\sB$-Baire-like and $\sB$ contains every pointwise convergent sequence, whose limit belong $\wG$, on open subgroups of $\wG$. Then every sequence $\{\varphi_n\}_{n \in \NN} \subseteq \widehat{G}$ that converges pointwise at every point of an open subgroup $H$ of $G$ is equicontinuous. As a consequence, every Baire-like group is $\omega$-barrelled.
\end{enumerate}
\ethm
\pf
(1) Let $\{B_n\}$ be an ng-$\sB$-barrel set in $G$. Then there is an open subgroup $H$ of $G$ such that $H \subseteq \bigcup_{n \in \NN} B_n$. We must prove that $\{B_n\}$ stabilizes.

Reasoning by contradiction, assume that $\{B_n\}$ does not stabilize. That is to say, for every $n \in \NN$ there is $m_n \in \NN$ such that $n<m_m$ and $B_n \subsetneq B_{m_n}$. Without loss of generality, we may assume that $m_n = n + 1$.

By hypothesis, for each $n \in \NN$, there is $K_n \in \sB$ such that $B_n = K_n^\lhd$ for all $n \in \NN$. Since $B_n \subsetneq B_{n+1}$, it follows that $K_{n} \nsubseteq K_{n+1}$. Indeed, if $K_{n} \subseteq K_{n+1}$ then
$B_{n+1}=K_{n+1}^\lhd\subseteq K_n^\lhd = B_n$, which is impossible.

Now, pick an arbitrary element $\varphi_n \in K_n \setminus K_{n+1}$ for each $n \in \NN$.
The sequence $\{\varphi_n\}$ converges to $0$ pointwise on $H$. Indeed, let $g \in H$ and $\epsilon > 0$ be arbitrary. Choose $m \in \NN$ such that $1/4m < \epsilon$. There exists $n_m \in \NN$ such that $jg \in B_n = K_n^\lhd$ for all $1 \leq j \leq m$ whenever $n \geq n_m$. Thus $\varphi_n(jg) \in [-\frac{1}{4}, \frac{1}{4}]$ for all $1 \leq j \leq m$ and $n \geq n_m$. This implies that $\varphi_n(g) \in [-\frac{1}{4m}, \frac{1}{4m}]$ for all $n \geq n_m$. Thus $\{\varphi_n(g)\}$ converges to $0$ for all $g \in H$. This yields a non-trivial convergent sequence on $H$, which is a contradiction that completes the proof.
\medskip

\noindent (2) Without loss of generality, we may assume that $\{\varphi_n\}_{n \in \NN}$ converges to $0$ on an open subgroup $H$. Set $K_n = \{\varphi_m\}_{m \geq n} \in \sB$ and $B_n := (\{\varphi_m\}_{m \geq n} \cup \{2\varphi_m\}_{m \geq n})^\lhd$. We first check that $\{B_n\}_{n \in \NN}$ is an ng-$\sB$-barrel set:
\begin{itemize}
    \item[(a)] Pick an element $x \in H$. Since $\{\varphi_n(x)\}$ converges to $0$, there is $n_x \in \NN$ such that $|\varphi_n(x)| \leq \frac{1}{8}$ for all $n\geq n_x$.
    Therefore, Lemma \ref{lem1} yields\,  $\varphi_n(2x)=2\varphi_n(x)$, which implies that $|(2\varphi_n)(x)| \leq \frac{1}{4}$ for all $n \geq n_x$.
    Thus $H \subseteq \bigcup_{n \in \NN} B_n$.
    \item[(b)] This inclusion is trivial.
\end{itemize}
Since $G$ is $\sB$-Baire-like, there is $n_0 \in \NN$ such that $\emptyset \neq \int (B_{n_0}) \subseteq \int (B_{n_0+1}) \subseteq \dots$. Now, pick $g \in B_{n_0}$ and a neighborhood $V$ of $0$ in $G$ such that $g + V \subseteq \int B_{n_0}$. We have that $|\varphi(g+v)| \leq \frac{1}{4}$ and $|(2\varphi)(g+v)| \leq \frac{1}{4}$ for all $\varphi \in K_{n_0}$ and $v \in V$. Therefore, $|\varphi(g+v)| \leq \frac{1}{8}$ for all $\varphi \in K_{n_0}$. In particular, $|\varphi(g)| \leq \frac{1}{8}$ for all $\varphi \in K_{n_0}$.

Now, for all $v \in V$, we have $|\varphi(g)| \leq \frac{1}{8}$ and $|\varphi(g+v)| \leq \frac{1}{8}$. Applying Lemma \ref{lem1}, we obtain $|\varphi(v)| \leq |\varphi(g+v)| + |\varphi(-g)| \leq \frac{1}{4}$. That is, $\varphi(v) \in [-\frac{1}{4}, \frac{1}{4}]$ for all $\varphi \in K_{n_0}$. This implies that $V \subseteq K_{n_0}^\lhd$, and consequently, $K_{n_0}^\lhd$ is a neighborhood of $0$ in $G$. This plainly implies that the sequence $\{\varphi_n\}_{n \in \NN}$ is an equicontinuous subset of $\widehat{G}$.
\epf
\medskip

We will need the following technical result in the sequel.

\bthm\label{convergencia} (convergencia)
Let $G$ be a topological abelian group and $H \le G$ an open subgroup. If $(\chi_n)_{n=1}^\infty$ is a sequence of continuous characters on $H$ converging pointwise to $0$ on $H$, then there exists a sequence of continuous characters $(\psi_n)_{n=1}^\infty$ on $G$ such that $\psi_n\vert{}_H = \chi_n$ for all $n \in \mathbb{N}$, and $(\psi_n)$ converges pointwise to $0$ on $G$.
\ethm
\pf\ Let $M$ be an intermediate subgroup such that $H \le M < G$, and let $(\psi_{n, M})_{n=1}^\infty \subset \widehat{M}$ be a sequence of continuous characters extending $\chi_n$ such that $\psi_{n, M}(m) \to 0$ for all $m \in M$.
For any element $g \in G \setminus M$, we extend the sequence to the strictly larger subgroup $M' = \langle M, g \rangle$ as follows:
If there is no $k\in\ZZ$ \st $kg\in M$, define $\psi_{n, M'}: M' \to \TT$ by
$\psi_{n, M'}(m + k \cdot g) = \psi_{n, M}(m)$ {for } $m \in M, k \in \mathbb{Z}$. Clearly, $\psi_{n, M'}$ is well defined.
Furthermore, for any $x = m + k \cdot g \in M'$, we have that
$\lim_{n \to \infty} \psi_{n, M'}(x) = \lim_{n \to \infty} \psi_{n, M}(m) = 0$.

Suppose now that for some $m\in\ZZ$ it holds that $mg\in M$ and let
$m_0 \ge 2$ be the smallest positive integer such that $m_0 \cdot g \in M$.
Choose the element $c_n\in\left(-\frac{1}{2m_0}, \frac{1}{2m_0}\right]$ \st $m_0c_n=\psi_{n, M}(m_0g)$ and
define $\psi_{n, M'}: M' \to \TT$ by
$\psi_{n, M'}(m + k \cdot g) = \psi_{n, M}(m) + \frac{k}{m_0} c_n$ {for } $m \in M, 0 \le k < m_0$,
which is clearly a group homomorphism. Furthermore, taking into account that $H$ is an open subgroup of $G$,
it is readily seen that $\psi_{n, M'}$ is continuous. On the other hand,
$\lim_{n \to \infty} \psi_{n, M'}(x) = \lim_{n \to \infty} \psi_{n, M}(m) + \frac{k}{m_0} \lim_{n \to \infty} c_n = 0$.
Hence, the sequence $(\psi_{n, M'})$ converges pointwise to $0$ on $M'$.

Now, apply Zorn's lemma to the set $\sS$ of all pairs $(M, (\psi_{n, M}))$, where $M$ is an intermediate subgroup $H \le M \le G$
and $(\psi_{n, M})$ extends $(\chi_n)$ while converging pointwise to $0$ on $M$, and the partial order defined by
$(M, (\psi_{n, M})) \le (M', (\psi_{n, M'}))$ if $M\leq M'$ and $\psi_{n, M'}$ extends $\psi_{n, M}$ for each $n$.
It is clear that $\sS$ is inductive, since for any chain of such subgroups and character sequences,
it suffices to take the union of the subgroups and defining the characters piecewise gives a well-defined upper bound.
Therefore, there is a maximal element $(L, (\psi_{n, L}))$ in $\sS$. The immediate conclusion is that $L$ must equal $G$
because, if an element $g \notin L$ existed, we could repeat the same arguments used above in order to obtain
a new element $(L', (\psi_{n, L'}))$ strictly larger than $(L, (\psi_{n, L}))$, contradicting maximality. \epf
\mkp

Theorems \ref{prp1.1} and \ref{convergencia} yield the following implications.   

\bcor\label{cor01}
Let $G$ be a topological abelian group and let $\sB$ be a bornology on $\wG$. If the dual group $\widehat{G}_p$ contains no non-trivial convergent sequences, 
then every ng-$\sB$-barrel set stabilizes and, as a consequence, the group $G$ is strongly $\sB$-Baire-like.
\ecor
\pf
Reasoning as in Theorem \ref{prp1.1} ({\rm 1}), if $\{B_n\}$ is an ng-$\sB$-barrel set in $G$ that does not stabilize, it follows that there is a non trivial sequence $\{\varphi_n\} \subseteq \widehat{G}$ that converges pointwise 
to $0$ on an open subgroup $H$ of $G$. Applying Theorem \ref{convergencia}, we may assume that the sequence $\{\varphi_n\}$ that converges pointwise
to $0$ on $G$. This contradicts the hypothesis and completes the proof.
\epf

\bcor\label{prp3}
Let $G$ be a MAP topological abelian group and let $\sB$ be a bornology on $\wG$ that contains all convergent sequences in $\wG_p$. Then $G$ is strongly $\sB$-Baire-like if and only if its dual $\widehat{G}_p$ contains no non-trivial convergent sequences.
\ecor
\pf By Corollary \ref{cor01}, we only need to verify that if $\widehat{G}_p$ contains a non-trivial convergent sequence, then $G$ contains an ng-$\sB$-barrel set that does not stabilize. Let $\{\varphi_n\}_{n \in \NN} \subseteq \widehat{G}$ be a non-trivial sequence that converges pointwise to $0$ on $G$. Defining $B_n = (\{\varphi_m\}_{m \geq n})^\lhd$, we obtain an ng-$\sB$-barrel set that does not stabilize.
\epf

\bthm\label{cor1} (cor1)
Let $G$ be a totally bounded abelian group and let $\sB$ be a bornology on $\wG$ that contains all convergent sequences in $\wG_p$. The following properties are equivalent:
\begin{itemize}
\item[(i)] $G$ is $\sB$-Baire-like.
\item[(ii)] The dual group $\widehat{G}_p$ contains no non-trivial convergent sequences.
\item[(iii)] $G$ is strongly $\sB$-Baire-like.
\item[(iv)] $G$ is dually Baire-like.
\end{itemize}
\ethm
\pf (i) $\Rightarrow$ (ii). If $G$ is $\sB$-Baire-like, then every convergent sequence in $\widehat{G}_p$ is equicontinuous and, since $G$ is totally bounded, finite.
\medskip

\noindent (ii) $\Rightarrow$ (iii). It suffices to apply Corollary \ref{cor01}. 

\noindent (iii) $\Rightarrow$ (i) is obvious.

\noindent (iv) $\Rightarrow$ (ii). Assume there is a non-trivial convergent sequence $\{\varphi_n\}$ in $\wG_p$. Defining $B_n = (\{\varphi_m\}_{m \geq n})^\lhd$, we obtain an ng-$\sB$-barrel set that does not stabilize. Furthermore, $(B_1 \cup \dots \cup B_n)^\rhd \supseteq \{\varphi_m\}_{m \geq n+1} \neq \{0\}$, which implies that $G$ is not dually Baire-like.

\noindent (i) $\Rightarrow$ (iv). Assume that $G$ is not dually Baire-like. Then there is a sequence $\{F_n\}$ of $\sK$-neighborhoods of $0_G$ such that $\bigcup_{n \in \NN} F_n$ is $G_\delta$-dense in $G$ and $(F_1 \cup \dots \cup F_n)^\rhd \neq \{0\}$ for all $n \in \NN$. Thus, there is a continuous character $0 \neq \varphi_1 \in F_1^\rhd$. Suppose we have already defined $m$ natural numbers $\{1 = n_1, \dots, n_m\}$ and $m$ non trivial continuous characters $\{\varphi_1, \dots, \varphi_m\} \subseteq (F_1 \cup \dots \cup F_{n_m})^\rhd$. Since $\varphi_m \neq 0$, there is an element $g \in G$ such that $|\varphi_m(g)| > 1/4$. Furthermore, since $\bigcup_{n \in \NN} F_n$ is $G_\delta$-dense in $G$, there is $x \in \bigcup_{n \in \NN} F_n$ such that $|\varphi_m(x)| = |\varphi_m(g)| > 1/4$. Let $n_{m+1}$ be the smallest natural number such that $x \in F_1 \cup \dots \cup F_{n_{m+1}}$ and select $0 \neq \varphi_{m+1} \in (F_1 \cup \dots \cup F_{n_{m+1}})^\rhd$. Note that the character $\varphi_{m+1} \neq \varphi_j$ for $1 \leq j \leq m$. Thus we obtain a non-trivial sequence $\{\varphi_m\}$. In order to prove that the sequence $\{\varphi_m\}$ converges to $0$ in $\wG_p$, it suffices to show that $G = \bigcup_{n \in \NN} (\{\varphi_m\}_{m \geq n})^\lhd$. Pick an arbitrary element $g \in G$. Because $\bigcup_{n \in \NN} F_n$ is $G_\delta$-dense in $G$, there is an element $x \in \bigcup_{n \in \NN} F_n$ such that $\varphi_m(x) = \varphi_m(g)$ for all $m \in \NN$. Take a natural number $n_m$ such that $x \in F_1 \cup \dots \cup F_{n_m} \subseteq (\{\varphi_m\}_{m \geq n_m})^\lhd$. It follows that $g \in (\{\varphi_m\}_{m \geq n_m})^\lhd$. This completes the proof.
\epf
\bigskip

It is clear that every Baire topological abelian group is $\sB$-Baire-like for every bornology $\sB \subseteq \wG$. However, there are examples of $\sB$-Baire-like groups that are not Baire. An example of this is presented below.

\bexm\label{exm1}
Let $G = C_p(\beta \NN, \TT)$ and let $\sB$ be a bornology on $\wG$ that contains all convergent sequences in $\wG_p$. Since $\widehat{G}_p = A_p(\beta\NN)$, which contains no non-trivial convergent sequences, by Theorem \ref{cor1}, the group $G$ is strongly $\sB$-Baire-like. However, the group $\widehat{G}_p$ contains infinite compact subsets, which means that $G_p$ is not Baire (see \cite{FerHerSepTri}).
\eexm
\section{$\sB$-Baire-like groups}
Here, we study some specific properties of $\sB$-Baire-like groups.

\bdfn\label{def2} (def2)
Let $G$ be a topological abelian group and let $\{A_n\}\subseteq 2^{\widehat{G}}$ denote a sequence of subsets.
It is said that $\{A_n\}$ converges pointwise to $0$ at every point of a subgroup $H$ of $G$ if
for every $\epsilon>0$ and $h\in H$ there is $n_0\in\NN$ \st $|\chi(h)|<\epsilon$ for all
$\chi\in A_n$ and $n\geq n_0$.
\edfn
\medskip

A variation of the arguments made in the proof of Theorem \ref{prp1.1} yields the following characterization of $\sB$-Baire-like groups.

\bthm\label{prp2} (prp2)
Let $G$ be a topological abelian group and let $\sB$ be a bornology in $\wG$. The group $G$ is $\sB$-Baire-like if and only if for every decreasing sequence of subsets $\{K_n\}\subseteq\sB$
that converges pointwise to $0$ at every point of an open subgroup of $G$ there is $n\in\NN$ such that $K_n$ is equicontinuous.
\ethm
\pf
\noindent \emph{Sufficiency:} Let $\{A_n\}$ be a ng-$\sB$-barrel set in $G$ whose union contains an open subgroup $H$.
Hence, for every $n\in\NN$, there is $K_n\in\sB$ \st $A_n=K_n^{\, \lhd}$.
We claim that the sequence of subsets $\{A_n^{\, \rhd}=K_n^{\, \lhd\rhd}\}$ converges pointwise to $0$ at every point of $H$.
Indeed, by hypothesis, we have that $A_n =A_n^{\, \rhd\lhd}$ for all $n\in\NN$.
Take an arbitrary element $g\in H$, $\epsilon >0$, and pick $m\in \NN$ \st $1/4m<\epsilon$.
There is $n_m\in\NN$ such that $jg\in A_n$ for all $1\leq j\leq m$ and $n\geq n_m$.
Thus $\varphi(jg) \in [-\frac{1}{4}, \frac{1}{4}]$ for all $\varphi\in A_n^{\, \rhd}$ whenever $n\geq n_m$ and $1\leq j\leq m$.
This implies that $\varphi(g) \in [-\frac{1}{4m}, \frac{1}{4m}]$ for all $\varphi\in A_n^{\, \rhd}$ and $n\geq n_m$.
Thus the sequence of subsets $\{A_n^{\, \rhd}\}$ converges pointwise to $0$ at every point $g \in H$.
Therefore, there is $n\in\NN$ \st $A_n^{\, \rhd}$ is equicontinuous. Hence $A_n= A_n^{\, \rhd\lhd}$ is a neighbourhood of $0$ in $G$, which implies that $\operatorname{int}(A_n)\neq\emptyset$.

\noindent \emph{Necessity:} Let $\{K_n\}\subseteq \sB$ be a decreasing sequence of subsets that converges pointwise to $0$ at every point of an open subgroup $H\subseteq G$.
Clearly, the sequence of subsets $\{K_n^{\, \lhd}\}$ is a ng-$\sB$-barrel set in $G$ whose union contains $H$.
Since the group $G$ is $\sB$-Baire-like, there is $n_0 \in \mathbb N$ \st $\operatorname{int}(K_{n_0}^{\, \lhd}) \neq \emptyset$.
Thus, there are $g \in G$ and $U \in \mathcal{N}_G(0)$ such that $g+U \subseteq \operatorname{int}(K_{n_0}^{\, \lhd})$. We also assume without loss of generality
that $jg\in K_{n_0}^{\, \lhd}$ for $1\leq j\leq 3$, which yields $|\varphi(g)|\leq \frac{1}{12}$.
Set $V\in \mathcal{N}_G(0)$ \st $V+V\subseteq U$. We have $\varphi(jv)=\varphi(g+jv-g)=\varphi(g+jv)-\varphi(g)$ for $1\leq j\leq 2$. As a consequence,
$\varphi(jv) \pmod{\mathbb Z} = \varphi(g+jv)-\varphi(g) \in [-\frac{1}{4}-\frac{1}{12}, \frac{1}{4}+\frac{1}{12}]$,
which yields $|\varphi(v)| \leq \frac{1}{3}$ and  $|\varphi(2v)| \leq \frac{1}{3}$. Under these conditions, it follows that $\varphi(2v)=2\varphi(v)$, which implies
$|\varphi(v)| \leq \frac{1}{6}$ for all $\varphi \in K_{n_0}$. That is, we have verified that $V\subseteq K_{n_0}^{\, \lhd}$, which
implies that $K_{n_0}$ is equicontinuous.
\epf
\medskip

\bqtn\label{qtn1} (qtn1)
When is it possible to replace sequences of subsets in $\sB$ by sequences of characters in Proposition \ref{prp2}?
That is, when is it true that a topological abelian group $G$ is $\sB$-Baire-like if and only if
every sequence $\{\varphi_n\}\subseteq \widehat{G}$ converging pointwise to $0$ at each point of an open subgroup of $G$ is equicontinuous?
\eqtn
\medskip

Below, we provide an affirmative answer to Question \ref{qtn1} for metrizable groups.

\bthm\label{thm2} (thm2)
Let $G$ be a metrizable group and let $\sB$ be a bornology on $\wG$ containing all pointwise convergent sequences. Then $G$ is $\sB$-Baire-like if and only if every pointwise convergent sequence on an open subgroup of $\widehat{G}$ is equicontinuous.
\ethm
\pf \noindent \emph{Necessity:} This fact has been proven in Theorem \ref{prp1.1}.

\noindent \emph{Sufficiency:} Let $\{U_n\}$ be a countable base of locally quasi-convex neighbourhoods of $0$ \st
$U_{n+1}\subseteq U_n$ for all $n\in\NN$ and
let $\{B_n\}$ be an arbitrary ng-$\sB$-barrel set in $G$ \st $\cup B_n$ contains an open subgroup $H$ of $G$. If there is $n_0 \in \mathbb N$ \st $B_{n_0}^{\vartriangleright}$ is equicontinuous,
then $B_{n_0}^{\vartriangleright \vartriangleleft} \in \mathcal{N}_{G}(0)$ and, since $B_{n_0} = B_{n_0}^{\vartriangleright \vartriangleleft}$, it follows that $B_{n_0} \in \mathcal{N}_{G}(0)$ and, as a consequence, $\operatorname{int}(B_{n_0})\neq \emptyset$.

Therefore, we may assume without loss of generality that $B_n^{\vartriangleright}$ is equicontinuous for no $n \in \mathbb N$.
Thus $B_n = B_n^{\vartriangleright \vartriangleleft} \not\supseteq U_n$ for all $n\in\NN$, which implies that $B_n^{\vartriangleright} \not\subseteq U_n^{\vartriangleright}$.

Take $\varphi_n \in B_n^{\vartriangleright} \setminus U_n^{\vartriangleright} \neq \emptyset$ for all $n \in \mathbb N$.
We claim that the sequence $\{\varphi_n\}$ converges pointwise to $0$ on $H$.
Indeed, take an arbitrary element $g \in H$. Since $H\subseteq \cup B_n$,  there is $n_0$ \st if $n \geq n_0$
then $g \in B_n$. We have $\varphi_n \in B_n^{\vartriangleright}$, which yields $|\varphi_n(g)| \leq \frac{1}{4}$ for all $n \geq n_0$.
This verifies that $\{\varphi_n\}$ converges pointwise to $0$ on $H$.
Hence, the subset
$\{\varphi_n\}$ is equicontinuous, which implies that $\{\varphi_n\}^{\vartriangleleft} \in \mathcal{N}_{G}(0)$.
Thus, there is $U_m \subseteq \{\varphi_n\}^{\vartriangleleft}$ such that $|\varphi_n(u)| \leq \frac{1}{4}$ for all $u \in U_m$ and $n \in \mathbb N$. 
In particular, $|\varphi_m(u)| \leq \frac{1}{4}$ for all $u \in U_m$ or, equivalently, $\varphi_m \in U_m^{\vartriangleleft}$.
This is a contradiction that completes the proof.
\epf

\section{$\go$-barrelled groups}

In this section we study some surprising properties of $\go$-barrelled groups. First, we need the following technical but simple lemma, whose proof we include for the reader's sake.

\blem\label{lem001} (lem001)
Let $K_1$ and $K_2$ be compact subsets of $\wG_p$ with $\we(K_1) \cdot \we(K_2) \leq \gk$. Then $\we(K_1 \pm K_2) \leq \gk$.
\elem
\pf Consider the \emph{addition/substraction} map $\phi\colon K_1\times K_2\to \wG_p$ and consider the dual map $\Phi\colon C_\infty(K_1\pm K_2)\to C_\infty(K_1\times K_2)$, which is a continuous embedding with respect to the supremum norm.
Since $\we(K_1\times K_2)\leq \gk$, it follows that the Banach space $C_\infty(K_1\times K_2)$ has a dense subset with cardinality less than or equal to $\gk$. Therefore, identifying $C_\infty(K_1\pm K_2)$ with its image in the metric space $C_\infty(K_1\times K_2)$, it follows that the former space also contains a dense subset of cardinality less than or equal to $\gk$. From this fact, it follows that $\we(K_1\pm K_2)\leq \gk$.
\epf
\medskip

From the above lemma, a characterization of $\gk$-barrelled groups follows. First, given a subset $A$ of an abelian group $G$,
set $A_{(n)} := \{g \in G : g, 2g,\dots ,ng \in A \}$ (see \cite{gabri:jmaa2023}).

\bthm\label{prp04} (prp04)
Let $G$ be a topological abelian group. Then the following assertions hold true:
\begin{enumerate}
 \item $G$ is g-barrelled iff every $\sK$-neighbourhood of $0$ has nonempty interior.
 \item $G$ is $\gk$-barrelled iff every $\sK$-neighbourhood $U$ of $0$ such that $U_{(2)}$ is $\gk$-dense has nonempty interior.
 \end{enumerate}
\ethm
\pf The proof of (1) is obvious. (2) \emph{Necessity}: Let $U$ be a $\gk$-dense $\sK$-neighbourhood of $0$. Then there is a compact subset $K\subseteq \wG_p$ \st $K^\triangleleft= U$. 
Furthermore, we may assume without loss of generality that $0\in K$.
We claim that $\we(K)\leq \gk$. Indeed, take the compact subset $K+K\subseteq \wG_p$ and $V=(K+K)^\triangleleft$. By \cite[Lemma 2.1]{gabri:jmaa2023}, we have that $V=U_{(2)}$. Hence,
there is a subset $A\subseteq G$ \st $|A|\leq \gl\leq \gk$ and $G=A+V$.
It will suffice to verify that $A$ separates the points in $K$. Now, if we take two elements $\chi_1\not=\chi_2$ in $K$, there is $g\in G$ \st $|\chi_1(g)-\chi_2(g)|>1/3$. Take $a\in A$ and $v\in V$ with $g=a+v$. Then we have that $|\chi_i(a)|\leq 1/2$ and $|\chi_i(v)|\leq 1/8$, $i\in\{1,2\}$.

Therefore, $1/3\leq |\chi_1(a)+\chi_1(v)-(\chi_2(a)+\chi_2(v))|\leq |\chi_1(a)-\chi_2(a)|+|\chi_1(v)-\chi_2(v)| \leq |\chi_1(a)-\chi_2(a)|+|\chi_1(v)|+|\chi_2(v)|\leq |\chi_1(a)-\chi_2(a)|+1/4$. Hence $1/12\leq |\chi_1(a)-\chi_2(a)|$, which implies that $A$ separates the points in $K$ and, as a consequence, the compact subset $K$ has weight less than or equal to $\gk$. Since $G$ is $\gk$-barrelled, it follows that $K$ is equicontinuous and, as a consequence, we have that $K^\triangleleft$ has nonempty interior.

\noindent (2) \emph{Sufficiency}: Let $K\subseteq \wG_p$ \st $\we(K)\leq \gl\leq \gk$. By Lemma \ref{lem001}, the metric space $C_\infty(K,\TT)$ has a dense subset $D$ with $|D|\leq \gl\leq \gk$. As a consequence, the metric space $G_{|K}\subseteq C_\infty(K,\TT)$ has a dense subset $A$ with $|A|\leq \gl\leq \gk$. From this fact, it follows that $K^\triangleleft$ is $\gk$-dense in $G$ and, therefore, a neighbourhood of $0$.\epf
\medskip

Next, we show how the class of $\omega$-barrelled groups coincides with the class of $\aleph_0$-barrelled groups; that is, the class of groups in which every metrizable compact subset of $\wG$ is equicontinuous. It is obvious that every $\aleph_0$-barrelled group is $\omega$-barrelled. Therefore, it will suffice to prove the converse implication. First we need the following result, whose proof is folklore but that we include for the reader's sake.

\blem\label{lem_metric}
Let $G$ be an abelian topological group and $K \subset G$ a metrizable compact subset. There exists an invariant pseudometric $d$ on $G$ such that the restriction $d|_K$ is a metric that induces the subspace topology on $K$.
\elem
\pf
Since $K$ is a metrizable compact space, the diagonal $\Delta = \{(x, x) : x \in K\}$ is a $G_\delta$ set in $K \times K$. Given the uniform structure of $G$, there exists a sequence $\{U_n\}_{n=0}^\infty$ of symmetric neighborhoods of $0 \in G$ satisfying:
\begin{enumerate}
    \item $U_{n+1} + U_{n+1} + U_{n+1} \subset U_n$ for all $n \in \mathbb{N}$,
    \item $\left( \bigcap_{n=1}^\infty U_n \right) \cap (K - K) = \{0\}$.
\end{enumerate}
Where $\bigcap_{n=1}^\infty U_n$
is a closed subgroup of $G$.

We define a scale function $f: G \to [0, \infty)$ by:
$$f(g) = \begin{cases} 2^{-n} & \text{if } g \in U_n \setminus U_{n+1}, \\ 1 & \text{if } g \notin U_0, \\ 0 & \text{if } g \in \bigcap_{n=1}^\infty U_n. \end{cases}$$

Following the Birkhoff-Kakutani construction, we define the invariant pseudometric $d: G \times G \to [0, \infty)$ as:
$$d(x, y) = \inf \left\{ \sum_{i=1}^k f(g_i) : x - y = \sum_{i=1}^k g_i, \text{ with } g_i \in G \right\}.$$
By construction, $d$ is invariant ($d(x+g, y+g) = d(x, y)$) and continuous on $G$. To see that $d|_K$ is a metric, note that if $x, y \in K$ with $x \neq y$, then $x - y \in (K - K) \setminus \{0\}$. By condition (2), $x - y \notin \bigcap_{n=1}^\infty U_n$, which implies $d(x, y) > 0$. Finally, the identity map $i: (K, \tau_G) \to (K, d)$ is a continuous bijection from a compact space to a Hausdorff space, hence a homeomorphism.
\epf
\medskip

We are now in position to prove our main result in this section.
\medskip

\noindent \textbf{Proof of Theorem \ref{Thm_A}:}
It will suffice to prove that for every metrizable compact subset $K$ of the topological abelian group $G$, there is a sequence $S:=\{w_n\}$ convergent to the neutral element of $G$ such that $S^\rhd \subseteq K^\rhd$. Thus, let $K$ be a metrizable compact subset of $G$. By Lemma \ref{lem_metric}, there is an invariant pseudometric $d$ of $G$ that defines the topology of $K$. Assume, without loss of generality, that $\operatorname{diam}(K) \le 1$ and divide $K$ into a finite number of closed balls of radius $1/4$ whose union contains $K$:
$$ K \subseteq \overline{B}(x_1, 1/4) \cup \dots \cup \overline{B}(x_{n(0)}, 1/4). $$
Take each of these balls $\overline{B}(x_j, 1/4)$ and divide it again into smaller closed balls of radius $\frac{1}{4^2}$:
$$ \overline{B}(x_j, 1/4) \subseteq \overline{B}(x_{j1}, \frac{1}{4^2}) \cup \dots \cup \overline{B}(x_{jn(j)}, \frac{1}{4^2}). $$
\textbf{Claim:} Observe that if $y \in \overline{B}(x, \frac{1}{4^m})$, then $d(x, y) \le \frac{1}{4^m}$. Therefore:
$$ d(0, 2^k(y-x)) \le d(0, y-x) + \dots + d(2^{k-1}(y-x), 2^k(y-x)) = 2^k d(0, y-x) \le \frac{2^k}{4^m}. $$
Reasoning by induction, let us suppose that we have defined the closed ball $\overline{B}(x_{t}, \frac{1}{4^m})$ for each finite chain $t = t_1 \dots t_m$ of natural numbers of length $l(t):=m$. Again, we divide it into a finite number of smaller closed balls:
$$ \overline{B}(x_{t}, \frac{1}{4^m}) \subseteq \overline{B}(x_{t1}, \frac{1}{4^{m+1}}) \cup \dots \cup \overline{B}(x_{tn(t)}, \frac{1}{4^{m+1}}). $$
Thus, we can define all closed balls $\overline{B}(x_{t}, \frac{1}{4^{m+1}})$ for chains of length $l(t)=m+1$, once we have defined the closed balls for all chains of length $m$. By the Claim above, if $y, x \in \overline{B}(z, \frac{1}{4^{m+1}})$, then:
$$ d(0, 2^k(y-x)) < \frac{2^k}{4^{m+1}}, \quad 1 \le k \le m+1. $$
Now we take the following sequence:
$$ S = \{x_1, 2x_1, \dots, x_{n(0)}, 2x_{n(0)}\} \cup \{x_t - x_{ti}, 2(x_t - x_{ti}), \dots, 2^{m+1}(x_t - x_{ti}) : l(t)=m, m \in \mathbb{N}\}. $$
We claim that this sequence converges to the neutral element. Indeed, take an arbitrary neighborhood $U$ of the neutral element. Since $K$ is compact, there exists $\delta > 0$ such that if $x, y \in K$ and $d(x, y) < \delta$, then $x - y \in U$. Pick a natural number $k_0$ such that $\frac{1}{2^{k_0}} < \delta$. For any chain $t$ with length $l(t) = k \ge k_0$, we have:
$$ d(0, 2^j(x_t - x_{ti})) \le \frac{2^j}{4^k} \le \frac{2^k}{4^k} = \frac{1}{2^k} \le \frac{1}{2^{k_0}} < \delta. $$
Therefore, $2^{j}(x_t - x_{ti}) \in U$ for $0 \le j \le l(t)$. This proves the convergence of $S$ to the neutral element of $G$.

Now, let $\varphi$ be a continuous character of $G$ such that $\varphi \in S^\rhd$. If we take an element $x_j$ in the first level, then $|\varphi(x_j)| \le \frac{1}{4}$ and $|\varphi(2x_j)| \le \frac{1}{4}$, which implies $|\varphi(x_j)| \le \frac{1}{2} \cdot \frac{1}{4} \le \frac{1}{4}$. Again, reasoning by induction, suppose we have verified that $|\varphi(x_{t})| \le \left( \frac{1}{2} + \dots + \frac{1}{2^{l(t)}} \right) \frac{1}{4} \le \frac{1}{4}$ for all chains of length $l(t)=m$. Then:
$$ |\varphi(x_{ti})| \le |\varphi(x_t)| + |\varphi(x_{ti} - x_t)| \le \left( \frac{1}{2} + \dots + \frac{1}{2^{l(t)}} \right) \frac{1}{4} + \frac{1}{2^{l(t)+1}} \frac{1}{4} \le \left( \frac{1}{2} + \dots + \frac{1}{2^{l(t)+1}} \right) \frac{1}{4} \le \frac{1}{4}. $$
Thus, $|\varphi(x_t)| \le \frac{1}{4}$ for all finite chains $t$. Since the set $\{x_t : 1 \le l(t) < \infty \}$ is dense in $K$, the continuity of $\varphi$ leads us to:
$$ |\varphi(x)| \le \frac{1}{4}, \quad \forall x \in K. $$
\epf
\bigskip

\noindent \textbf{Proof of Corollary A:}
If in the previous theorem we take a dense subset $D$ of $K$, there is no problem in replacing the centers of the balls $x_t$ with points located in the dense subset $D$. In such a case, we would possibly need to take a larger number of closed balls in each step. But the essential facts of the construction remain without further modification.
\medskip

Theorem A also provides a characterization of g-barrelled groups in the realm of separable metrizable groups.

\bcor\label{thm22} (thm22)
A separable, metrizable group $G$ is g-barrelled if and only if it is Baire-like.
\ecor
\pf \noindent \emph{Necessity:} If $G$ is g-barrelled then it is also $\go$-barrelled and we apply Theorem \ref{thm2}.

\noindent \emph{Sufficiency:} If $G$ is Baire-like then every convergent sequence in the dual group $\wG_p$ is equicontinuous and, applying Theorem A, every metrizable compact subset of $\wG_p$ is also equicontinuous. Now, since $G$ is separable, it follows that every compact subset in $\wG_p$ is metrizable and, therefore, equicontinuous.
\epf
\medskip

\noindent \textbf{Proof of Theorem \ref{Thm_C}:}
Let $\{\chi_n\}$ be a Cauchy sequence in $\wG_\sK$. Clearly the sequence $\{k\chi_n\}_{n<\omega}$ is also Cauchy for every $k<\omega$. Since $G$ is $\sigma$-compact, there is an increasing sequence of compact subsets $\{K_m\}_{m<\go}$ \st $G=\bigcup_{m<\go}K_m$.
Hence, given any $m<\omega$, there is $n(m)<\omega$ such that
\begin{equation*}
\left\| k(\chi_{n(m)}-\chi_{n(m)+p}) \right\| _{K_{m}}:=\sup \{|k(\chi_{n(m)}(x)-\chi_{n(m)+p}(x))|: x \in K_{m}\}<1/2^{2m}
\end{equation*}
for all $p<\go$ and $1\leq k\leq 2^{m+1}$.

We now define $\psi_{m}:=\chi_{n(m)}$ for all $m<\omega$ and set
$$S:=\{\psi_{1},2\psi_{1}\} \cup \bigcup_{m<\omega} \bigcup_{k=1}^{2^{m+1}} \{k(\psi_{m+1}-\psi_{m})\}.$$

Since $\left\| k(\psi_{m+1}-\psi_{m}) \right\| _{K_{m}}<1/2^{2m}$ for all $1\leq k\leq 2^{m+1}$, it follows that $S$ converges pointwise to $0$. We now apply ideas similar to the proof of Theorem \ref{Thm_A}.
Take any element $g$ in $S^\lhd$. Since $\psi_{1}$ and $2\psi_{1}$ belong to $S$, it follows that $|(2\psi_1)(g)|=|2\psi_1(g)|\leq 1/4$ and, as a consequence, $|\psi_1(g)|\leq 1/(2\cdot 4)$.
Analogously, we have
$$\{\psi_{m+1}-\psi_{m},\dots,k(\psi_{m+1}-\psi_{m}),\dots,2^{m+1}(\psi_{m+1}-\psi_{m})\}\subset S,$$
so that, using an inductive argument, we obtain
$$|\psi_{m+1}(g)-\psi_{m}(g)|\leq 1/(2^{m+1}\cdot 4)$$
for all $m<\omega$. Then we have that
$$\psi_{m+1}=\psi_{1}+(\psi_{2}-\psi_{1})+\dots+(\psi_{m+1}-\psi_{m})$$
with $|\psi_{1}(g)|\leq 1/(2\cdot 4)$ and $|\psi_{k+1}(g)-\psi_{k}(g)|\leq 1/(2^{k+1}\cdot 4)$ for $1\leq k\leq m$.
Hence
$$|\psi_{m+1}(g)|\leq |\psi_{1}(g)|+\sum_{k=1}^{m}|\psi_{k+1}(g)-\psi_{k}(g)|\leq \frac{1}{4}\sum_{k=1}^{m}\frac{1}{2^{k+1}}<\frac{1}{4}.$$
This means that $\psi_{m+1}\in S^{\lhd\rhd}$ for all $m<\omega$. On the other hand, the initial Cauchy sequence $\{\chi_n\}$ must converge pointwise to some element $\chi_0\in \operatorname{Hom}(G,\TT)$. Thus, the subsequence $\{\psi_{m}\}=\{\chi_{n(m)}\}$ also converges pointwise to $\chi_0$. This means that $|\chi_0(g)|\leq \frac{1}{4}$ for all $g\in S^\lhd$.

Now, since $G$ is $\go$-barrelled, $S$ is an equicontinuous subset of $\wG$, which yields the continuity of $\chi_0$. Thus, the Cauchy sequence $\{\chi_{n}\}$ converges pointwise to the continuous character $\chi_0$, which clearly implies convergence in $\wG_\sK$. Finally, it is a well-known fact that if $G$ is hemicompact, then $\widehat{G}_\sK$ is metrizable.
\epf

\section{An $\go$-barrelled metrizable group that is not $g$-barrelled}\label{Javier's question}

et $X = \omega_\omega$ be the first limit ordinal of uncountable cardinality, which has countable cofinality. For each $n < \omega$, we define the closed ordinal intervals
\[
K_n = [\omega_n + 1, \omega_{n+1}] \subset \omega_\omega,
\]
which partition $\omega_\omega$ into clopen, compact blocks.

\bdfn\label{compactification w_w}
Let $C_b(X)$ be the $C^*$-algebra of all bounded complex-valued continuous functions on $X$, and let $\mathcal{A}$ be the subalgebra of functions $f \in C_b(X)$ satisfying the following condition: there exists $n_f < \omega$ such that $f(x) = f(\omega_n)$ for all $x \in K_n$ and all $n > n_f$.

The subalgebra $\sA$ canonically yields the following natural embedding:
\begin{equation*}
\ggm : X \hookrightarrow \prod\limits_{f\in \sA}\overline{f(X)} \quad \text{with} \quad \ggm(x) = (f(x))_{f\in \sA}
\end{equation*}
\mkp
This embedding defines a compactification of $X$ denoted by $(\ggm X, \ggm)$, where $\ggm X := \overline{\ggm(X)}$ and
every $f\in\sA$ can be extended to acontinuous fuction defined on $\ggm X$. We denote this continous extension also by $f$ for simplicity's sake.
\edfn

Next, we identify the compactification remainder $\ggm X \setminus X$.

\medskip

\blem\label{Remainder} (Remainder)
The remainder $K = \gamma X \setminus X$ is homeomorphic to $\beta\omega \setminus \omega$.
\elem

\pf
Let $Y := \{\omega_n : n < \omega\}$. We first prove that $\ggm X \setminus X = \text{cl}_{\ggm X} Y \setminus Y$.

Indeed, let $p \in \ggm X \setminus X$ be an arbitrary point in the remainder of $\ggm X$. For every finite subset $F = \{f_1, \dots, f_l\} \subseteq \sA$ and $\epsilon > 0$, the set
\[
V(p,F, \epsilon) = \left\{ x \in \gamma X : |f_j(p) - f_j(x)| < \epsilon, \ 1 \leq j \leq l \right\}
\]
is a neighborhood of $p$, and the family $\{ V(p,F,\epsilon) : F \in \sA^{<\omega}, \epsilon > 0\}$ forms a neighborhood base at $p$. Since $p \in \ggm X \setminus X$, it follows that for every neighborhood $V(p,F,\epsilon)$, there is an infinite sequence $(n_m) \subseteq \omega$ such that $A_m := V(p,F,\epsilon) \cap K_{n_m} \neq \emptyset$. Set $B := \{\omega_{n_m} : m \in \omega\}$. By the definition of $\sA$, there exists a natural number $n_F < \omega$ such that $f_j(x) = f_j(\omega_n)$ for all $x \in K_n$ whenever $n > n_F$. Therefore, $V(p,F,\epsilon) \cap B \neq \emptyset$, which means that $p$ is also an accumulation point of $Y$. The converse inclusion is obvious.

Thus, in order to verify that $\ggm X \setminus X$ is homeomorphic to $\beta\omega \setminus \omega$, it suffices to prove that $\text{cl}_{\ggm X} Y$ is a canonical copy of $\beta \omega$. We show this by proving that every complex-valued bounded function $f \colon Y \to \CC$ can be extended to a continuous function $\bar{f} \colon \ggm X \to \CC$.
To do this, we define the map $h \colon X \to \omega$ by setting $h(K_n) = \{n\}$ and $g \colon \omega \to \CC$ by $g(n) = f(\omega_n)$ for all $n < \omega$. Note that $g \circ h \colon X \to \CC$ satisfies $(g \circ h)_{|Y} = f$. Furthermore, the map $g \circ h$ belongs to $\sA$ since it is constant in each block $K_n$ (with $n_{g \circ h} = 0$). Therefore, it can be extended to a continuous map $\bar{f} \colon \ggm X \to \CC$. Since $\bar{f}_{|Y} = f$, this completes the proof.
\epf
\medskip

We now construct the following metrizable topological group.

\bexm[The Metrizable Topological Group $G$]
Let $G$ be the subgroup of $C(\gamma X, \mathbb{T})$ consisting of all functions $f \in \mathcal{A}$ that are constant on each block $K_n$ for all $n < \omega$. We equip $G$ with the topological group topology $\tau$ generated by uniform convergence on the countable family of compact subsets:
\[
\mathcal{F} = \{ K_n \}_{n < \omega} \cup \{ K \},\ \text{where}\ K=\ggm X\setminus X.
\]
Because $\mathcal{F}$ is a countable family of compact sets whose union covers the domain $\gamma X$, the topology $\tau$ is Hausdorff and is generated by a countable family of sup-pseudometrics. Thus, $(G, \tau)$ is a metrizable topological group.
\eexm
\medskip

We can now prove the main features of the group $G$

\bthm
$G$ is a $\go$-barrelled metrizable group that is not $g$-barrelled.
\ethm
\pf
It will suffice to prove the following facts:

Let $\widehat{\ggm X} = \{ \delta_x \in \wG \mid x \in \ggm X \}$ be the set of evaluation characters in the dual group corresponding to points in $\ggm X$, equipped with the pointwise convergence topology $\tau_p$. Then
\begin{enumerate}
    \item $G$ is $\go$-barrelled.
    \item $\widehat{\cl_{\ggm X}Y}$ is $\tau_p$-compact.
    \item $\widehat{cl_{\ggm X} Y}$ is not equicontinuous on $G$.
\end{enumerate}
\medskip

\noindent (1). First, observe that
the topology of $G$ is stronger than the topology it inherits from the pointwise convergence topology of the group $C_p(\gamma X, \mathbb{T})$ but weaker than the one it inherits from the uniform topology on $C_c(\gamma X, \mathbb{T})$. Since $\gamma X$ is a zero-dimensional compact space, it follows that all these groups have the same algebraic dual $A(\ggm X)$, the free topological abelian group generated by $\ggm X$ (see \cite{pestov95,GalHer99,Aussenhofer}). On the other hand, since every $f \in G$ is uniquely determined by its values on $Y = \{\omega_n : n < \omega\} \approx \omega$, $G$ is topologically isomorphic to $C(\beta\omega, \mathbb{T})$. Thus, its dual group is isomorphic to 
$A(\beta\omega)$ and also isomorphic to $<\cl_{\ggm X}Y>$, which is a retract subgroup of $A(\ggm X)$. Let $(w_n)$ be a sequence in $A(\beta\omega)$ weakly converging to $0$. Because the free topological abelian group $A(\beta\omega)$ strongly respects compactness (see \cite{GalHer:fundamenta}), it follows that $(w_n)$ converges to $0$ in the natural (Markov) topology of $A(\beta\omega)$. However, the group $A(\ggm X)$ contains no non trivial convergent sequences, which implies that $G$ is $\go$-barrelled.
\medskip

\noindent (2). The evaluation mapping $\delta: \gamma X \to \widehat{G}_p$ defined by $p \mapsto \delta_p$ is continuous by the definition of the product topology on $\mathbb{T}^G$. Since $\cl_{\ggm X}Y$ is compact, its image  $\widehat{\cl_{\ggm X}Y} = \delta(\cl_{\ggm X}Y)$ is $\tau_p$-compact.
\medskip

\noindent (3). For $\widehat{\cl_{\ggm X}Y}$ to be equicontinuous on $G$, the topology $\tau$ would have to coincide with the uniform topology on the entire space $\cl_{\ggm X}Y$. However, $\tau$ is strictly coarser. We can construct a sequence of functions $(f_n)_{n<\omega} \subset G$ that converges to $0$ uniformly on every individual $K_m$ and on $K$, but maintains $\sup_{y \in \gamma X} |f_n(y) - 0| = 1/2$.

Indeed, identifying $\mathbb{T}$ with $\mathbb{R}/\mathbb{Z}$, define $f_n \colon \ggm X \to \mathbb{T}$ by setting $f_n(x) = 1/2 \pmod 1$ if $x \in K_n$, and $f_n(x) = 0 \pmod 1$ if $x \notin K_n$ (which includes all other blocks $K_m$ for $m \neq n$ and the remainder $K$).

Clearly, for any fixed compact set in our family $\mathcal{F}$ (either an individual block $K_m$ or the remainder $K$), the sequence $(f_n)$ is eventually identically zero. Thus, $(f_n)$ converges to $0$ in the topology $\tau$ of $G$. However, because $\sup_{y \in \gamma X} |f_n(y) - 0| = 1/2$ for all $n$, the sequence does not converge uniformly on the whole space $\gamma X$. This shows that $\widehat{\ggm X}$ is not equicontinuous on $G$, completing the proof.
\epf

\section{Unordered Baire-like groups}
Our first result shows that the property of being unordered Baire-like is preserved by quotients.

\blem\label{lem0001} (lem0001)
Let $G$ be an unordered Baire-like group and let $\phi\colon G\to H$ be a quotient group homomorphism.
Then $H$ is also unordered Baire-like.
\elem
\pf Let $\{F_n\}$ be a sequence of $\sK$-neighbourhoods of $0_H$ with empty interior in $H$.
Since $\phi$ is a quotient, it follows that $\phi^{-1}(F_n)$ has empty interior for all $n\in\NN$.
Thus $\bigcup_{n\in\NN}\phi^{-1}(F_n)\neq G$, which yields $\bigcup_{n\in\NN}F_n \neq H$. \epf
\bigskip

\blem\label{lema0001.1} (lema0001.1)
Let $G$ be a totally bounded abelian group and let $\widehat{G}_p$ denote its dual group equipped with the finite-open topology.
If $\widehat{G}_p$ is unordered Baire-like, then $G$ contains no nontrivial convergent sequence.
As a consequence, every compact subset of $G$ is either finite or uncountable.
\elem
\pf Apply Corollary \ref{cor1} since every unordered Baire-like group is Baire-like. \epf
\bigskip

A continuous surjection $f:L\to C$ is said to be {\em irreducible} if whenever $A\subseteq L$ is closed and proper, then $f[A]\neq C$. The proof of the following result is a variation of the arguments used in \cite[Lemma 2.8]{FerHerSepTri}. Our formulation of this well-known result requires that the compact subsets are placed in a topological group and contain the neutral element, but the same proof can easily be adapted to fulfill these requirements.

\blem\label{lem0002} (lem0002)
Let $K$ be a compact subset of an Abelian topological group $G$ that is not scattered and contains the neutral element as an accumulation point.
Then $K$ contains a compact subset $L$, containing the neutral element as an accumulation point, \st there exists a continuous irreducible surjection $f : L \to \TT$, satisfying $f(0)=0$, and there is a countable family $\{V_n\}$ of open subsets of $L$, each containing a compact subspace $L_n$ with $|L_n|\geq \cc$, \st for every open set $V$ of $L$ there is $n\in\NN$ with $L_n\subseteq V_n \subseteq V$. Furthermore, in case $K$ is $0$-dimensional, we can replace $\TT$ by $\{0,1\}^\omega$.
\elem
\pf Set $C=\TT$ for a general compact non-scattered space $K$ and set $C=\{0,1\}^\omega$ in case we are only concerned with a $0$-dimensional compact space. Proving that there is a closed subspace $L$ of $K$ possessing an irreducible continuous surjection $f:L\to C$ is folklore. Nevertheless, we include a sketch of the proof for the reader's sake.

Indeed, fix a continuous surjective map $f:K \to C$ satisfying $f(0)=0$ (which exists by \cite[Lemma 2.7]{FerHerSepTri}, composing with a translation in $\TT$ if necessary) and let
\[\sF:=\{X \subseteq K: 0\in X,\ f(0)=0,\ X \mbox{ is closed and } f[X]=C\}.\]
Plainly, this set is inductively ordered and the maximal element has the desired properties (see also \cite[Prob. 3.1.C(a)]{engel}).

Let $\{B_n\}$ be a countable base of $C$. We have that each $B_n$ contains a compact subspace, say $K_n$, with $|K_n|=\cc$ for all $n \in \NN$. Let $V_n:=f^{-1}[B_n]$ and $L_n:=f^{-1}[K_n]$. Then $L_n$ is compact, $|L_n|\geq\cc$, and $L_n\subseteq V_n$ for all $n \in \NN$.

Let $V$ be an open subset of $L$. Then $L\setminus V$ is closed in $L$, hence $F:=f[L\setminus V]$ is closed in $C$. Since $f$ is irreducible, $B:=C \setminus F$ is non-empty, and obviously open, hence there is a $B_n$ contained in $B$. It follows that $V_n \subseteq f^{-1}[B] \subseteq V$, as required.
\epf
\medskip

\blem\label{lema0006.1} (lema0006.1)
Let $G$ be a totally bounded abelian group \st every uncountable compact subset of $G$ contains an uncountable independent subset. If the dual group $\widehat{G}_p$, equipped with the finite-open topology, is unordered Baire-like, then every compact subset of $G$ is finite.
\elem
\pf Let $K$ be an infinite compact subspace of $G$, which we assume contains the neutral element as an accumulation point without loss of generality. If $K$ were scattered, it would contain a convergent sequence, which is impossible (see \cite[Th. 4]{MRS}). We can therefore assume that $K$ is not scattered. Applying Lemma \ref{lem0002}, $K$ contains a compact subset $L$, containing the neutral element as an accumulation point, \st there exists a continuous irreducible surjection $f : L \to \TT$, satisfying $f(0)=0$, and there is a countable family $\{V_n\}$ of open subsets of $L$, each containing a compact subspace $L_n$ with $|L_n|\geq \cc$, \st for every open set $V$ of $L$ there is $n\in\NN$ with $L_n\subseteq V_n \subseteq V$.

Set $F_n:=L_n^\rhd$ for all $n\in\NN$. According to the hypotheses we have assumed, each $L_n$ contains an uncountable independent subset, say $A_n$. As a consequence, it is easily seen that $F_n$ has empty interior for all $n\in\NN$ (see \cite[Lemma 2.6]{FerHerSepTri}).

Since $\widehat{G}_p$ is unordered Baire-like, it follows that $\bigcup_{n\in\NN} F_n\subsetneq \wG$. Now, pick $\phi \in \wG\setminus \bigcup_{n\in\NN} F_n$ and consider $U=\{g \in G: |\phi(g)|<1/4\}$. Then there is $n\in\NN$ \st $L_n\subseteq U$. However $\phi\notin F_n$, which implies there exists an element $w\in L_n$ \st $|\phi(w)|>1/4$. This is a contradiction that completes the proof.
\epf
\bigskip

We can now prove our main result in this section. It characterizes totally bounded abelian groups containing no nontrivial compact subset.
\medskip

\noindent \textbf{Proof of Theorem \ref{Thm_D}:}
\emph{Necessity:} Since every compact subset of $G$ is finite, every $\sK$-neighbourhood of $0_G$ is actually a neighbourhood of $0_G$.

\noindent \emph{Sufficiency:} This is proved analogously to the proof of \cite[Theorem 2.12]{FerHerSepTri}.
\epf
\bigskip

The next result provides a characterization of g-barrelled totally bounded abelian groups.

\bcor\label{cor0001} (cor0001)
Let $G$ be a totally bounded abelian group. The following assertions are equivalent:
\begin{enumerate}
  \item[(i)] $G$ is g-barrelled.
  \item[(ii)] $G$ is unordered Baire-like.
\end{enumerate}
\ecor
\pf First observe that, since $G$ is totally bounded, $G$ is equipped with the pointwise convergence topology on $\wG_p$. Furthermore, $G\cong \widehat{(\wG_p)}_p$ (see \cite{RaczTrig}). Therefore, we can exchange the roles of the groups $G$ and $\wG$ in Theorem D.

\noindent (i)$\Rightarrow$ (ii). Since $G$ is g-barrelled and totally bounded, all compact subgroups of $\wG$ are finite, which means that every $\sK$-neighbourhood of $0_G$ is actually a neighbourhood of $0_G$ and, therefore, has nonempty interior.

\noindent (ii)$\Rightarrow$ (i). If $G$ is unordered Baire-like then $\wG_p$ contains no infinite compact subset. As a consequence, $G$ is g-barrelled.
\epf
\section{Acknowledgement}

We are very grateful to Professors Lydia Au$\beta$enhofer, Xabier Dom\'inguez and Javier Trigos-Arrieta for sharing helpful comments with us.


\begin{thebibliography}{99}


\bibitem{Aussenhofer}
{ L. Au{\rm \ss}enhofer,} {\em Contributions to the Duality Theory
of Abelian Topological Groups and to the Theory of Nuclear
Groups,} Dissertation. T\"ubingen 1998; Dissertationes
    Mathematicae (Rozprawy Matematyczne) CCCLXXXIV, Polska Akademia
    Nauk, Instytut Matematyczny, Warszawa, 1999.
	
\bibitem{AusDik:2021}
Au\ss enhofer, L. and Dikranjan, D., \emph{Mackey groups and Mackey topologies}, Dissertationes Mathematicae 567 (2021), 1--141.

\bibitem{ausdikgio}
Au\ss enhofer, L. and Dikranjan, D., Giordano Bruno,  A., \emph{Topological groups and the Pontryagin–van Kampen duality—an introduction},
De Gruyter Stud. Math., 83, De Gruyter, Berlin, (2022).

\bibitem{Banaszczyk}  \textsc{W. Banaszczyk}, \textit{Additive subgroups of topological
vector spaces}, \textrm{Lecture Notes in Math}. {\bf 1466}. Springer--Verlag,
Berlin-Heidelberg. 1991.


\bibitem{chasco}
\noindent M. J. Chasco, {\em Pontryagin duality for metrizable groups},
Arch. Math. { 70} (1998), 22--28 .
	
\bibitem{ChascoDominguezTkachenko}
Chasco, M. J., Dom\'\i nguez, X. and Tkachenko, M., {\em Duality properties of bounded torsion topological Abelian groups,}
J. Math. Anal. Appl. 448 (2017), 968-981.

\bibitem{ChMPVa:99}
Chasco, M. J., Mart\'in Peinador, E. and Tarieladze, V., \emph{On Mackey Topology for groups},
Stud. Math. 132, No.3, (1999), 257--284.	


\bibitem{comractri:04}
\noindent W. W. Comfort and S. U. Raczkowski and F. J.
Trigos-Arrieta, {\em The dual group of a dense subgroup},
Czechoslovak Math. Journal { 54}, No. 129, (2004), 509--533.
	
\bibitem{ComfortRoss1964}
Comfort, W. W. and Ross, K. A., {\em Topologies induced by groups of characters,} {Fundamenta Mathematicae}, {55} (1964), {283--291}


\bibitem{engel}
Engelking. R., {\em General Topology,} Heldermann Verlag, Berlin 1989.

\bibitem{FerHerTka} Ferrer, M. V., Hern\'andez, S. and Tkachenko, M., {\em On convergent sequences in dual groups},
RACSAM 114, No. 71, (2020).

\bibitem{FerHerSepTri}
Ferrer, M.V., Hernández, S., Sepúlveda, I., Trigos-Arrieta, F. J.,
{\em The Baire property and precompact duality,}
Proc. Amer. Math. Soc. Ser. B 11 (2024), 612--623.

\bibitem{FHU:jmaa}
M. V. Ferrer, S. Hern\'andez, and V. Uspenskij, \emph{Precompact groups
and property (T)}, J. Math. Anal. Appl. 404, No. 2, (2013) 221--230.

\bibitem{gabri:jmaa2023}
Gabriyelyan, S., {\it Compatible group topologies on a locally quasi-convex abelian group and the Mackey group problem},
J. Math. Anal. Appl. 524, No. 2, (2023), Paper No. 127111, 22 pp.


\bibitem{GalHer:fundamenta} Galindo, J. and Hern\'andez, S., {\em The concept of boundedness and the Bohr compactification of a MAP Abelian group}, Fund. Math. 159, No. 3, (1999) 195218.
	
\bibitem{GalHer99}
Galindo, J. and Hern\'andez, S., {\em Pontryagin-van Kampen reflexivity for free Abelian \tg s,} Forum Math. 11 (1999), 399--415.
	




	\bibitem{HR}
	Hewitt, E. and Ross K. A., {\em Abstract Harmonic Analysis: Volume I, Structure of Topological Groups,
		Integration Theory, Group Representations,} Grundlehren der mathematischen Wissenschaften,
		Springer New York, 1963.

\bibitem{hofmor_compactgroups}
Hofmann, K. H., Morris, S. A., \emph{The structure of compact groups—a primer for the student—a handbook for the expert},
De Gruyter Stud. Math., 25, Fifth Edition, De Gruyter, Berlin, (2023).

\bibitem{mackey}
{G. Mackey}, \textit{Harmonic analysis as the explotation of
symmetry: a historical survey,} {\textrm Bull. Amer. Math. Soc.} { 3} (1980),
543--698.

\bibitem{morr}  {S. A. Morris}, \textit{Pontryagin Duality and the
Structure of Locally Compact Abelian groups}, Cambridge University Press.
1977.

\bibitem{MRS}
Mrowka, S., Rajagopalan, M. and Soundararajan, T., {\em A characterization
of compact scattered spaces through chain limits}, TOPO 72 -- General
Topology and its Applications, Lecture Notes in Mathematics 378 (1974), 288--297.

\bibitem{Noble:tams}
Noble, N., {\em $k$-Groups and Duality},Trans. Amer. Math. Soc.,  151, No. 2,
(1970), 551--561

\bibitem{pestov95}  \noindent V. Pestov, \emph{Free Abelian topological groups
and the Pontryagin-van Kampen duality}, Bull. Austral. Math. Soc. {52}, No. 2, (2009), 297--311.

\bibitem{RaczTrig}
Raczkowski, S. U. and Trigos-Arrieta, F. J., \emph{Duality of totally bounded Abelian groups},
Bolet\'\i n de la Sociedad Matem\'atica Mexicana, 3, No. 7, (2001), 1--12.

\bibitem{rudin}  {W. Rudin}, \textit{Fourier {A}nalysis on
{G}roups}. Interscience,
New York. 1962.

\bibitem{Saxon:MA1972}
Saxon, S., \emph{Nuclear and product spaces, Baire-like spaces, and the strongest locally
convex topology}, Math. Ann. 197 (1972), 87--106.

\bibitem{Trigos:2025}
F. J: Trigos-Arrieta, \emph{Sequentially Barrelled Groups}, Algebra, Topology, and Dynamical Systems,
A conference in honor of Dikran Dikranjan on the occasion of his 75th birthday and of Ivan Prodanov on the occasion of the 90th anniversary of his birth.


\bibitem{Vilenkin}
{Vilenkin, N. Ya.}, \emph{The theory of characters of topological {A}belian groups with
              boundedness given}, {Izvestiya Akad. Nauk SSSR. Ser. Mat.},
 {15} (1951), {439--462}.
\end{thebibliography}
\end{document}